\documentclass[11pt]{article}

\usepackage{amsmath, amsthm, amssymb}
\usepackage{enumerate}
\usepackage{fullpage}

\usepackage{tikz}

\theoremstyle{plain}
\newtheorem{thm}{Theorem}[section]
  \theoremstyle{plain}
  \newtheorem{prop}[thm]{Proposition}
  \theoremstyle{definition}
  \newtheorem{defn}[thm]{Definition}
  \theoremstyle{plain}
  \newtheorem{lem}[thm]{Lemma}

\newcommand{\BBE}{\mathbb{E}}
\newcommand{\BFP}{\mathbb{P}}

\renewcommand{\baselinestretch}{1.05}

\author{
Choongbum Lee \thanks{Department of Mathematics, UCLA, Los
Angeles, CA, 90095. Email: choongbum.lee@gmail.com. Research
supported in part by a Samsung Scholarship.}
\and
Benny Sudakov \thanks{Department of Mathematics, UCLA, Los Angeles, CA 90095.
Email: bsudakov@math.ucla.edu.
Research supported in part by NSF grant DMS-1101185,
NSF CAREER award DMS-0812005 and by USA-Israeli BSF grant.}
}
\date{}

\begin{document}

\title{Dirac's theorem for random graphs}

\maketitle

\begin{abstract}
A classical theorem of Dirac from 1952 asserts that every graph
on $n$ vertices with minimum degree at least $\lceil n/2 \rceil$
is Hamiltonian. In this paper we extend this result to random graphs. Motivated by the study
of resilience of random graph properties
we prove that if $p \gg \log n /n$, then a.a.s.~every subgraph
of $G(n,p)$ with minimum degree at least $(1/2+o(1))np$
is Hamiltonian. Our result improves on previously known
bounds, and answers an open problem of
Sudakov and Vu.
Both, the range of edge probability $p$
and the value of the constant $1/2$ are asymptotically best possible.
\end{abstract}

\section{Introduction}

A \emph{Hamilton cycle} of a graph is a cycle which passes through
every vertex of the graph exactly once, and a graph is \emph{Hamiltonian}
if it contains a Hamilton cycle. Hamiltonicity is one of the most
central notions in graph theory, and has been intensively studied by
numerous researchers. The problem of determining Hamiltonicity
of a graph is one of the NP-complete problems that Karp listed
in his seminal paper \cite{Karp}, and accordingly,
one cannot hope for a simple classification of such graphs.
Therefore it is important to find general sufficient conditions for
Hamiltonicity and in the last 60 years many interesting results were obtained
in this direction. One of the first results of this type is a
classical theorem proved by
Dirac \cite{Dirac} in 1952, which asserts that every graph
on $n$ vertices of minimum degree at least $\lceil n/2 \rceil$
is Hamiltonian.

In this paper, we study Hamiltonicity of random graphs.
The model of random graphs we study is the binomial model $G(n,p)$
(also known as the Erd\H{o}s-Renyi random graph),
which denotes the probability space whose points are graphs with
vertex set $[n] = \{1,\ldots,n\}$ where each pair of vertices forms
an edge randomly and independently with probability $p$. We say that
$G(n,p)$ possesses a graph property $\mathcal{P}$
\textit{asymptotically almost surely}, or a.a.s. for brevity, if the
probability that $G(n,p)$ possesses $\mathcal{P}$ tends to 1 as $n$
goes to infinity. The earlier results on Hamiltonicity of random
graphs were proved by P\'osa \cite{Posa}, and Korshunov \cite{Korshunov}.
Improving on these results, Bollob\'as \cite{Bollobas}, and Koml\'os and Szemer\'edi \cite{KoSz} proved
that if $p\ge(\log n + \log \log n + \omega(n))/n$ for some function
$\omega(n)$ that goes to infinity together with $n$, then $G(n,p)$ is
a.a.s. Hamiltonian. The range of $p$ cannot be improved, since if $p \le (\log n + \log \log n - \omega(n))/n$,
then $G(n,p)$ a.a.s.~has a vertex of degree at most one.

Recently, in \cite{SuVu} the authors proposed to study Hamiltonicity
of random graphs in more depth by measuring how strongly the random
graphs possess this property. Let $\mathcal{P}$ be a monotone
increasing graph property. Define the \emph{local resilience} of a
graph $G$ with respect to $\mathcal{P}$ as the minimum number $r$
such that by deleting at most $r$ edges from each vertex of $G$, one
can obtain a graph not having $\mathcal{P}$. Using this notion, one
can state the aforementioned Dirac's theorem as ``$K_n$ has local
resilience $\left\lfloor n/2 \right\rfloor$ with respect to
Hamiltonicity''. Sudakov and Vu \cite{SuVu} initiated the systematic
study of resilience of random and pseudorandom graphs with respect
to various properties, one of which is Hamiltonicity. In particular,
they proved that if $p > \log^4 n / n$, then $G(n,p)$ a.a.s.~has
local resilience at least $(1/2+o(1))np$ with respect to
Hamiltonicity.

There are several other papers that studied local resilience of
random graphs with respect to various properties. For example,
Balogh, Csaba, and Samotij \cite{BaCsSa} studied the property of
containing almost spanning trees of bounded degree, and
B{\"o}ttcher, Kohayakawa, and Taraz \cite{BoKoTa} studied the
property of containing almost spanning subgraphs of bounded degree.
There is also a similar concept called {\em global resilience},
where one measures the total number of edges that needs to be
removed in order to obtain a graph without the given property.
Global resilience of random graphs has also been studied for many
properties; in fact, some results in this direction were obtained
even before the concept has been first formalized in \cite{SuVu}.
Haxell, Kohayakawa, and {\L}uczak \cite{HaKoLu1, HaKoLu2} studied it
with respect to the property of containing a fixed length cycle,
Dellamonica, Kohayakawa, Marciniszyn, and Steger \cite{DeKoMaSt}
studied it with respect to the property of containing long cycles,
and Alon and Sudakov \cite{AlSu} studied it with respect to
increasing the chromatic number. Recently, Conlon and Gowers
\cite{CoGo}, and Schacht \cite{Schacht} independently obtained a
breakthrough result which resolves several open problems in this
area, one of which establishes the global resilience of random
graphs with respect to containing fixed subgraphs.
For other recent results on resilience, see \cite{AlSu, BaLeSa, BKS,
HaLeSu, KiSuVu, KrLeSu, LeSa}.

The above mentioned result of Sudakov and Vu can be viewed as a generalization
of Dirac's Theorem, since a complete graph is also a random graph
$G(n,p)$ with $p=1$. This connection is very natural and in fact
most of the resilience results can be viewed as a generalization of
some classic graph theory result to random and pseudorandom graphs.
Note that, the constant $1/2$ in the resilience bound for Hamiltonicity cannot be further improved.
To see this, consider a partition of the vertex set of a random graph into two parts of size $n/2$
and remove all the edges between these parts.  Since the graph is random this removes roughly half of the edges
incident with each vertex and makes the graph disconnected.
However, things become unclear when one considers the range of $p$.
Recall that Bollob\'as, and Koml\'os and Szemer\'edi's result
mentioned above implies that if $p > C\log n/n$ for some $C>1$, then
$G(n,p)$ is a.a.s. Hamiltonian. Therefore it is
natural to believe, as was conjectured in \cite{SuVu}, that $G(n,p)$ has local resilience $(1/2+o(1))np$
with respect to Hamiltonicity already when $p \gg \log n/n$.

In addition to \cite{SuVu},
several other results have been obtained on this problem.
Frieze and Krivelevich \cite{FrKr} proved that there exist
constants $C$ and $\varepsilon$
such that if $p \ge C\log n/n$, then $G(n,p)$
a.a.s. has local resilience at least $\varepsilon np$
with respect to Hamiltonicity. This result was improved
by Ben-Shimon, Krivelevich, and Sudakov \cite{BeKrSu1} who showed
that for all $\varepsilon$, there exists a constant $C$,
such that if $p \ge C\log n/n$, then $G(n,p)$ has
local resilience at least $(1/6 -\varepsilon)np$. In
their recent paper \cite{BeKrSu2}, the same authors further improved this bound
to $(1/3-\varepsilon)np$.
Our main theorem completely solves the resilience problem of Sudakov and Vu.

\begin{thm}
\label{thm:mainthm}For every positive $\varepsilon$, there exists
a constant $C=C(\varepsilon)$ such that for $p\ge\frac{C\log n}{n}$,
a.a.s. every subgraph of $G(n,p)$ with minimum degree at least $(1/2+\varepsilon)np$
is Hamiltonian.
\end{thm}
As mentioned above, the constant $1/2$ and the range of edge probability $p$ are
both asymptotically best possible.

\bigskip

\noindent \textbf{Notation}. A graph $G=(V,E)$ is given by a pair of
its vertex set $V=V(G)$ and edge set $E=E(G)$. We sometimes use
$|G|$ to denote the order of the graph. For a subset $X$ of
vertices, we use $e(X)$ to denote the number of edges within $X$,
and for two sets $X,Y$, we use $e(X,Y)$ to denote the number of
edges $\{x,y\}$ such that $x \in X, y\in Y$ (note that
$e(X,X)=2e(X)$). We use $N(X)$ to denote the collection of vertices
of $V\setminus X$ which are adjacent to some vertex of $X$. For two
graphs $G_1$ and $G_2$ over the same vertex set $V$, we define their
intersection as $G_1 \cap G_2 = (V, E(G_1)\cap E(G_2))$, their union
as $G_1 \cup G_2 = (V, E(G_1) \cup E(G_2))$, and their difference as
$G_1 \setminus G_2 = (V, E(G_1) \setminus E(G_2))$

When there are several graphs under consideration, to avoid
ambiguity, we use subscripts such as $N_{G}(X)$ to indicate the
graph that we are currently interested in. We also use subscripts
with asymptotic notations to indicate dependency. For example,
$\Omega_\varepsilon$ will be used to indicate that the hidden
constant depends on $\varepsilon$. To simplify the presentation, we
often omit floor and ceiling signs whenever these are not crucial
and make no attempts to optimize absolute constants involved. We
also assume that the order $n$ of all graphs tends to infinity and
therefore is sufficiently large whenever necessary. All logarithms
will be in base $e \approx 2.718$.

\section{Properties of random graphs}

In this section we develop some properties of random graphs. The
following concentration result, Chernoff's bound (see, e.g., \cite[Theorem A.1.12]{AlSp}),
will be used to establish these properties.
\begin{thm} Let $\varepsilon$ be a positive
constant. If $X$ be a binomial random variable
with parameters $n$ and $p$, then
\[ \BFP\big(|X - np| \geq \varepsilon np\big) \leq e^{-\Omega_\varepsilon(np)}. \]
Also, for $\lambda \ge 3np$,
\[ \BFP\big(X - np \geq \lambda \big) \leq e^{-\Omega(\lambda)}. \]
\end{thm}

We also state another useful concentration result that we
will use later (see, e.g., \cite[Theorem 2.10]{JaLuRu}).
Let $A$ and $A'$ be sets such that $A' \subseteq A$.
Let $B$ be a fixed size subset of $A$ chosen
uniformly at random. Then the distribution of the random variable
$|B \cap A'|$ is called the \textit{hypergeometric distribution}.
\begin{thm} \label{lemma_hypergeometric}
Let $\varepsilon$ be a fixed positive constant and
let $X$ be a random variable with hypergeometric distribution. Then,
\[ \BFP\big(|X - \mathbb{E}[X]| \geq \varepsilon \mathbb{E}[X]\big) \leq e^{-\Omega_\varepsilon(\mathbb{E}[X])}. \]
\end{thm}

We first state two standard results on random graphs,
which estimates the number of edges and the degree
of vertices. We omit their proofs which consist of
straightforward applications
of Chernoff's inequality.
\begin{prop}
\label{prop:degreeandedges}For every positive $\varepsilon$, there exists
a constant $C$ such that for $p\ge\frac{C\log n}{n}$, the random
graph $G=G(n,p)$ a.a.s. has $e(G)=(1+o(1))\frac{n^2p}{2}$ edges, and
$\forall v \in V,\, (1-\varepsilon)np\le\deg(v)\le(1+\varepsilon)np$.
\end{prop}

\begin{prop}
\label{prop:edges} Let $p \ge \log n/ n$, and
$\omega(n)$ be an arbitrary function
which goes to infinity as $n$ goes to infinity.
Then in $G=G(n,p)$, a.a.s. for every
two subsets of vertices $X$ and $Y$,
\[ e(X,Y) = |X||Y|p + o(|X||Y|p + \omega(n) n). \]
\end{prop}

It is well-known that random graphs have certain
expansion properties, and that these properties are very useful in
proving Hamiltonicity. Next proposition shows that
the expansion property still holds even after removing some of
its edges. Similar lemmas appeared in \cite{KrLeSu, SuVu}.
\begin{prop}
\label{prop:expansion_smallset}For every positive $\varepsilon$,
there exists a constant $C$ such that for
$p\ge \frac{C\log n}{n}$, the random graph $G=G(n,p)$ a.a.s. has the following property.
For every graph $H$ of maximum degree at most $(\frac{1}{2} - 2\varepsilon)np$,
the graph $G' = G-H$ satisfies the following:
%Every subgraph $G'=G-H$, such that the
%maximum degree of $H$ is at most $(\frac{1}{2} - \varepsilon)np$, satisfies that:
\begin{enumerate}[(i)]
  \setlength{\itemsep}{1pt} \setlength{\parskip}{0pt}
  \setlength{\parsep}{0pt}
\item $\forall X \subseteq V\,,\,|X| \le (\log n)^{-1/4}p^{-1},\quad|N_{G'}(X)|\ge \left(\frac{1}{2}+\varepsilon\right)|X|np$,
\item $\forall X \subseteq V\,,\, n(\log n)^{-1/2} \le |X| \le \frac{\varepsilon }{2}n,\quad|N_{G'}(X)|\ge\left(\frac{1}{2}+\varepsilon\right)n$, and
\item $G'$ is connected.
\end{enumerate}
\end{prop}
\begin{proof}
Let $H$ be a graph of maximum degree at most $(\frac{1}{2} -
2\varepsilon)np$, and let $G'=G-H$.

\noindent $(i)$ To prove $(i)$, it suffices to prove that a.a.s. for all $X \subseteq V$
of size at most $(\log n)^{-1/4}p^{-1}$,
\[ |N_G(X)| \ge \left(1-\varepsilon\right)|X|np, \]
since it will imply by the maximum degree condition of $H$ that
\[ |N_{G'}(X)| \ge |N_G(X)| - \left(\frac{1}{2} - 2\varepsilon\right)np \cdot |X|
\ge \left(\frac{1}{2} + \varepsilon \right)np\cdot |X|. \]

Fix a set $X \subseteq V$ of size $|X| \leq (\log n)^{-1/4}p^{-1}$.
For each $v \in V \setminus X$, let $Y_v$ be indicator random
variable of the event that $v \in N(X)$. We have $\BFP(Y_v = 1) = 1
- (1-p)^{|X|} = (1+o(1))|X|p$ (the estimate follows from the fact
$|X|p=o(1)$). Let $Y=|N(X)|=\sum_{v \in V \setminus X} Y_v$ and note
that
\[ \mathbb{E}[Y] = \sum_{v \in V \setminus X} \BFP(Y_v =1) = (n-|X|)(1+o(1))|X|p = (1+o(1))|X|np. \]
Since the events $Y_v$ are mutually independent, we can apply the Chernoff's inequality to get
$\BFP\big( |Y - \mathbb{E}[Y]| \geq (\varepsilon/2) \mathbb{E}[Y]\big) \leq e^{-\Omega_\varepsilon(\mathbb{E}[Y])}$.
Combine this with the estimate on $\mathbb{E}[Y]$ and we have,
\[ \BFP\big( Y \le (1 - \varepsilon) |X|np  \big) \le e^{-\Omega_\varepsilon(\mathbb{E}[Y])} = e^{-\Omega_\varepsilon(|X|np)}, \]
for large enough $n$.
Since $np \ge C\log n$, the probability that $Y=|N(X)| < (1 - \varepsilon)|X|np$ is $e^{-\Omega_\varepsilon(|X|np)}=n^{-C'|X|}$,
where $C'=C'(\varepsilon, C)$ can be made arbitrarily large by choosing constant $C$ appropriately.

Taking the union bound over all choices of $X$, we get
\[
 \sum_{1 \leq |X| \leq (\log n)^{-1/4}p^{-1}} n^{-C'|X|} \le
 \sum_{k=1}^{n} \binom{n}{k} n^{-C'k}
 \leq \sum_{k=1}^{n} \left(\frac{en}{k} \cdot n^{-C'}\right)^k = o(1), \nonumber \]
which establishes our claim.

\bigskip

\noindent (ii) We will first prove that a.a.s. for every pair of disjoint sets
$X, Y$ of sizes $n (\log n)^{-1/2} \le |X| \le \frac{\varepsilon n}{2}$ and $|Y| \ge \left(\frac{1}{2}-\frac{3\varepsilon}{2}\right)n$,
\begin{equation}
\label{eq1}
e_G(X, Y) \ge \left(1 - \frac{\varepsilon}{2}\right)|X||Y|p
> \left(\frac{1}{2} - 2\varepsilon\right)|X|np.
\end{equation}
Indeed, let $X,Y$ be a fixed pair of disjoint sets
such that $n (\log n)^{-1/2} \le |X| \le \frac{\varepsilon n}{2}$ and $|Y| \ge \left(\frac{1}{2}-\frac{3\varepsilon}{2}\right)n$.
Then $\mathbb{E}[e_G(X,Y)] = |X||Y|p$ and by Chernoff's
inequality,
\[
\BFP\Big(e_G(X, Y) \le (1 - \varepsilon/2)|X||Y|p\Big) < e^{-\Omega_\varepsilon(|X||Y|p)} \le e^{-\Omega_\varepsilon(n (\log n)^{1/2})}.
\]
Since there are at most $2^{2n}$ possible choices of the pairs $X,
Y$ and the probability above is $\ll 2^{-2n}$, taking the union
bound will give our conclusion.

Condition on the event that (\ref{eq1}) holds, and assume that there
exists a set $X$ of size $n (\log n)^{-1/2} \le |X| \le
\frac{\varepsilon n}{2}$ which has less than
$(\frac{1}{2}+\varepsilon)n$ neighbors in $G'$. Then there exists a
set $Y$ of size at least $|Y| \ge n - (\frac{1}{2}+\varepsilon)n -
|X| \ge (\frac{1}{2}-\frac{3\varepsilon}{2})n$ disjoint from $X$
such that there are no edges between $X$ and $Y$ in $G'$. However,
this gives us a contradiction to (\ref{eq1}) since
\[ 0=e_{G'}(X,Y) \ge e_{G}(X,Y) - \left(\frac{1}{2} - 2\varepsilon\right)np\cdot |X| > 0. \]

\bigskip

\noindent (iii) Condition on the event that $(i)$ and $(ii)$ holds,
and assume that $G'$ is not connected. Let $X$ be a set of vertices
which induces a connected component in $G'$, and let $Y= V \setminus
X$. By part $(i)$, we know that $|X| \ge (\log n)^{-1/4}p^{-1}\cdot
\frac{np}{2}=\frac{1}{2}n(\log n)^{-1/4}$, and then by part $(ii)$,
we know that $|X|
> \frac{n}{2}$. On the other hand, since $Y$ must also contain a
connected component, the same estimate must hold for $Y$ as well.
However this cannot happen since the total number of vertices is
$n$. Therefore, $G'$ is connected.
\end{proof}

\section{Rotation and extension}

Our main tool in proving Hamiltonicity is
P\'osa's rotation-extension technique (see \cite{Posa} and
\cite[Ch. 10, Problem 20]{Lovasz}). We start by briefly
discussing this powerful tool
which exploits the expansion property of a graph, in order
to find long paths and/or cycles.

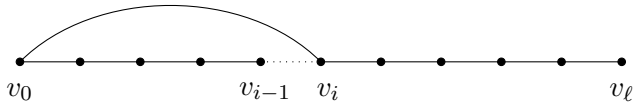
\begin{figure}[b]
  \centering
%  \begin{tabular}{c}
    \begin{tikzpicture}

%%
%% Second
%%
\def\xendtwo{8 cm}
\def\yshifttwo{2.2 cm}
%   \draw [color=white] (0,0) circle(0.1mm);

   \foreach \x in {0.0,0.8,...,8.8} {
            \draw [fill=black] (\x cm, \yshifttwo) circle (0.5mm);
        }

   \foreach \x in {0.0,0.8,1.6,2.4} {
            \draw (\x cm,\yshifttwo) -- (\x cm + 0.8cm,\yshifttwo);
        }
   \draw [dotted] (3.2 cm,\yshifttwo) -- (4.0cm,\yshifttwo);

   \foreach \x in {4.0,4.8,...,8.0} {
            \draw (\x cm,\yshifttwo) -- (\x cm + 0.8cm,\yshifttwo);
        }

 %          \draw [dashed] (1.5 cm,\yshift) -- (1.9cm,\yshift);
 %          \draw [dashed] (4.5 cm,\yshift) -- (4.9cm,\yshift);

    \draw (0.0cm, \yshifttwo) .. controls (1.0cm,\yshifttwo + 1.0cm)
            and (3.0cm,\yshifttwo + 1.0cm) .. (4.0cm, \yshifttwo);
%    \draw[-latex] (1.6cm, \yshifttwo) ..controls(3cm, \yshifttwo + 0.7cm)
%            and (4.2cm,\yshifttwo + 0.7cm) .. (5.6cm, \yshifttwo + 0.1cm);

    \draw (4.1cm, \yshifttwo - 0.4cm) node {$v_{i}$};
    \draw (3.26cm, \yshifttwo - 0.4cm) node {$v_{i-1}$};

%    \draw (5.6cm, \yshifttwo - 0.4cm) node {$v_{j}$};
%    \draw (4.8cm, \yshifttwo - 0.4cm) node {$v_{j-1}$};

    \draw (0.0cm, \yshifttwo - 0.4cm) node {$v_{0}$};
    \draw (8.0cm, \yshifttwo - 0.4cm) node {$v_{\ell}$};

\end{tikzpicture}
%  \end{tabular}
  \caption{Rotating a path}
  \label{fig_rotation}
\end{figure}

Let $G$ be a connected graph and let $P = (v_0, \cdots, v_\ell)$ be
a path on some subset of vertices of $G$ ($P$ is not necessarily a
subgraph of $G$). If $\{v_0, v_\ell\}$ is an edge of $G$, then we
can use it to close $P$ into a cycle. Since $G$ is connected, either
the graph $G \cup P$ is Hamiltonian, or there exists a longer path
in this graph. In the second case, we say that we \emph{extended}
the path $P$.

Assume that we cannot directly extend $P$ as above, and assume that
$G$ contains an edge of the form $\{v_0, v_i\}$ for some $i$. Then
$P' = (v_{i-1}, \dots, v_0, v_i, v_{i+1}, \dots, v_{\ell})$ forms
another path of length $\ell$ in $G \cup P$ (see Figure
\ref{fig_rotation}). We say that $P'$ is obtained from $P$ by a
\emph{rotation} with {\em fixed endpoint} $v_0$, \emph{pivot point}
$v_i$, and {\em broken edge} $\{v_{i-1}, v_i\}$. Note that after
performing this rotation, we can now close a cycle of length $\ell$
also using the edge $\{v_{i-1}, v_\ell\}$ if it exists in $G$. As we
perform more and more rotations, we will get more such candidate
edges (call them \emph{closing edges}). The rotation-extension
technique is employed by repeatedly rotating the path until one can
find a closing edge in the graph $G$ thereby extending the path.

Let $P''$ be a path obtained from $P$ by several rounds of
rotations. An important observation which we later will use is that
for an interval $I = (v_j, \dots, v_k)$ of vertices of $P$ ($0 \le
j < k \le \ell$), if no edges of $I$ were broken during these
rotations, then $I$ appears in $P''$ either exactly as it does in
$P$, or in the reversed order.

\medskip

We will use rotations and extension as described above to prove our
main theorem. The main technical twist is to split the given graph
into two graphs, where the first graph will be used to perform
rotations and the second graph to perform extensions. Similar ideas,
such as \emph{sprinkling}, has been used in proving many results on
Hamiltonicity of random graphs. The one which is closest to our
implementation, appears in the recent paper of Ben-Shimon,
Krivelevich, and Sudakov \cite{BeKrSu2}.

In the following two
subsections, we prove that our random graph indeed contains
subgraphs which can perform these two roles of rotation and extension.
All the graphs that we study from now on are defined
over the same vertex set, and we will use this fact
without further mentioning.

\subsection{Rotation}

\begin{defn}
\label{def:property_rotext}Let $\delta$ be a positive constant. We
say that a connected graph $G$ on $n$ vertices has property
$\mathcal{RE}(\delta)$ if one of the following holds for every path
$P$ (not necessarily a subgraph of $G$): (i) there exists a path
longer than $P$ in the graph $G\cup P$, or (ii) there exists a set
of vertices $S_{P}$ of size at least $|S_{P}|\ge \delta n$ such that
for every vertex $v\in S_{P}$, there exists a set $T_{v}$ of size
$|T_{v}|\ge \delta n$ such that for every $w\in T_{v}$, there exists
a path of the same length as $P$ in $G\cup P$ over the vertex set
$V(P)$ which starts at $v$ and ends at $w$.
\end{defn}

Informally, a graph has property $\mathcal{RE}(\delta)$ if every
path is extendable to a longer path, or can be rotated in many
different ways. The next lemma, which is the most crucial ingredient
of our proof, asserts that we can find a graph with property
$\mathcal{RE}(\frac{1}{2}+\varepsilon)$ in random graphs even after
deleting some of its edges.

\begin{lem}
\label{lem:rotation}For every positive $\varepsilon$, there exists
a constant $C=C(\varepsilon)$ such that for $p\ge\frac{C\log n}{n}$,
the random graph $G=G(n,p)$ a.a.s. has the following property:
for every graph $H$ of maximum degree at most $(\frac{1}{2} - 2\varepsilon)np$,
the graph $G'=G-H$ satisfies $\mathcal{RE}(\frac{1}{2} + \varepsilon)$.
%Every subgraph $G'=G-H$, such that the maximum degree of $H$ is at most $(1/2-2\varepsilon)np$
%satisfies $\mathcal{RE}(\frac{1}{2} + \varepsilon)$.
\end{lem}
\begin{proof}
Let $C$ be a sufficiently large constant such that the assertions of
Propositions \ref{prop:degreeandedges}, \ref{prop:edges}, and
\ref{prop:expansion_smallset} a.a.s.~hold, and
condition on all of these events. Let $H$
be a subgraph of $G(n,p)$ which has maximum degree at most
$(\frac{1}{2} - 2\varepsilon)np$, and let $G'=G-H$. By Proposition
\ref{prop:degreeandedges}, we know that $G'$ has minimum degree at
least $(\frac{1}{2}+\varepsilon)np$, and by Proposition
\ref{prop:expansion_smallset} (iii), we know that $G'$ is connected.
We want to show that $G'\in\mathcal{RE}(\frac{1}{2}+\varepsilon)$
for all choices of $H$. Consider a path $P=(v_{0},\dots,v_{\ell})$.
If there exists a path longer than $P$ in $G\cup P$, then there is
nothing to prove. Thus we may assume that this is not the case. For
a set $Z\subseteq V(P)$, let $Z^{+}={\{v_{i+1}|v_{i}\in Z\}}$ and
$Z^{-}={\{v_{i-1}|v_{i}\in Z\}}$. For a vertex $z$, let $z^-$ be the
vertex in $\{z\}^-$ and similarly define $z^+$.

\medskip
\noindent \textbf{Step 1 }: {\em Initial rotations.}
\smallskip

First we show that there exists a set $X$ of linear size such that
for all $v\in X$, there exists a path of length $\ell$ starting at
$v$ and ending at $v_{l}$. Such $X$ will be constructed iteratively.
In the beginning, let $X_{0}={\{v_0\}}$. Now suppose that we have
constructed sets $X_{i}$ of sizes $4^{-i}(np)^{i}$ up to some
nonnegative $i$. If $4^{-i}(np)^{i} \le \max\{1, (\log
n)^{-1/4}p^{-1}\}$, then either by the minimum degree of $G'$ (in
case, when $|X_i|=1$) or by Proposition
\ref{prop:expansion_smallset} (i), we know that $|N_{G'}(X_{i})|\ge
\left(\frac{1}{2}+\varepsilon\right)|X_{i}|np$. We must have
$N_{G'}(X_{i})\subseteq P$ as otherwise we can find a path longer than
$P$. Consequently, we can rotate the endpoints $X_{i}$ using the
vertices in $N_{G'}(X_{i})$ as pivot points. If a vertex $w \in
N_{G'}(X_{i})$ does not belong to any of $X_j$, $X_j^-$, $X_j^+$ for
$j < i$, then both edges of the path $P$ incident with $w$ were not
broken in the previous rotations. Hence, using $w$ as a pivot point,
we get either $w^-$ or $w^+$ as a new endpoint (see the discussion
at the beginning of the section). Therefore, at most two such pivot
points can give rise to a same new endpoint, and we obtain a set
$X_{i+1}$ of size at least
\begin{align*}
|X_{i+1}| &\ge \frac{1}{2} \left( |N_{G'}(X_{i})| - 3\sum_{j=0}^{i} |X_j| \right) \\
&\ge \frac{1}{2} \left( \left(\frac{1}{2} + \varepsilon \right) \left(\frac{np}{4}\right)^{i} np - o((np)^{i+1})\right)
\ge \left(\frac{1}{4}\right)^{i+1}(np)^{i+1},
\end{align*}
where $\sum_{j=0}^{i} |X_j| = o((np)^{i+1})$ since $|X_j| = (np/4)^j$ and $np \ge C\log n$.
Repeat the argument above until at step $t$ we have a set of endpoints $X_t$ of size at least $(\log n)^{-1/4}p^{-1}$, and
redefine $X_{t}$ by arbitrarily taking a subset of this set of size $|X_t|=\max\{1, (\log n)^{-1/4}p^{-1}\}$
(note that $t\le\frac{\log n}{\log(np/4)}\le\frac{\log n}{\log\log n}$
for $C\ge 4$). Apply the same argument as above to $X_t$ to find a
set of endpoints of size at least $\max\left\{\frac{np}{4}, \frac{n}{\log^{1/2} n} \right\}$.
Again, if necessary, redefine $X_{t+1}$ to be an arbitrary subset of this set of size
$|X_{t+1}|=n/\log^{1/2} n$, and repeat the argument above one more time, now using the
second part of Proposition \ref{prop:expansion_smallset} instead
of the first part to get $|N_{G'}(X_{t+1})| \ge (\frac{1}{2} + \varepsilon)n$. In the end, we obtain
a set $X_{t+2}$ of size at least
\[
|X_{t+2}| \ge \frac{1}{2}\left(|N_{G'}(X_{t+1})| - 3\sum_{j=1}^{t+1}|X_t|\right) \ge \frac{n}{4}.
\]

\medskip
\noindent \textbf{Step 2} : {\em Terminal rotation.}
\smallskip

Let $X = X_{t+2}$ be the set of size at least $\frac{n}{4}$ that
we constructed in Step 1.
We will show that another round of
rotation gives at least $(\frac{1}{2}+\varepsilon)n$ endpoints.
Let $Y$ be the set of all endpoints that we obtain by rotating $X$
one more time (note that $Y$ can contain vertices from $X$).

Partition the path $P$ into $k=\frac{\log n}{(\log\log n)^{1/2}}$
intervals $P_{1},\dots,P_{k}$ whose lengths are either
$\left\lfloor \frac{|P|}{k} \right\rfloor$ or $\left\lceil \frac{|P|}{k} \right\rceil$.
Every vertex $w\in X$ was obtained by $t+2$ rotations which
broke $t+2$ edges of $P$. If the interval $P_{i}$
contains none of these edges then the path from
$w$ to $v_l$ must traverse $P_{i}$
exactly in the same order as in $P$, or in the reverse order (see the discussion
at the beginning of the section). Let $\hat{X_i}$ be
the collection of vertices of $X$ which were obtained by rotation with some broken edges in
$P_{i}$. Let $X_{i,+}$ and $X_{i,-}$ be the vertices of $X$
such that $P_i$ is unbroken and paths from these vertices
to $v_l$ traverses $P_{i}$ in the original, and reverse order, respectively.
Note that $X=\hat{X_{i}}\cup X_{i,+}\cup X_{i,-}$ for all $i$.

The first key fact that we will now verify is that the set $\hat{X_i}$ is small for
most indices. Let $J$ be the collection of indices which
have $|\hat{X_i}| \ge (\log \log n)^{-1/4}|X|$. Since each vertex in
$X$ is obtained by at most $\frac{\log n}{\log \log n} + 2<2 \frac{\log n}{\log \log n}$
rotations, we can count the total number of broken edges used for constructing all
the points of $X$ in two ways to get
\[ |J|\cdot (\log \log n)^{-1/4}|X| \le |X| \cdot \frac{2\log n}{\log \log n}, \]
which implies $|J| \le 2\log n / (\log \log n)^{3/4} = o(k)$.

Our second key fact is that for a vertex $v_{j} \in P_i$
and a vertex $x \in X_{i,+}$, if $\{x, v_{j+1}\}$ is an edge of $G'$,
then $v_j \in Y$ (similarly, for $x \in X_{i,-}$, if $\{x, v_{j-1}\}$
is an edge of $G'$, then $v_j \in Y$). Therefore,
for all $i$, there are no edges of $G'$ between
$X_{i,+}$ and $(P_i \cap P_i^+) \setminus Y^+$, and between
$X_{i,-}$ and $(P_i \cap P_i^- )\setminus Y^-$. We will show that
if $|Y| < (\frac{1}{2} + \varepsilon)n$, then this cannot happen
because we will have to remove too many edges incident to $X$ from the
graph $G$ to form $G'$.

The number of edges incident to $X$ that we need to remove is at least,
\[
e_{G}(X, V\setminus P)+\sum_{i=1}^{k}\Big(e_{G}(X_{i,-}, (P_i \cap P_i^- )\setminus Y^-)+e_{G}(X_{i,+}, (P_i \cap P_i^+) \setminus Y^+)\Big).\]
Since $(P_i \cap P_i^- )\setminus Y^-$ and $(P_i \setminus Y)^-$  differs by at most one element (similar for $P_i^+$),
the above expression is
\[
e_{G}(X, V\setminus P)+\sum_{i=1}^{k}\Big(e_{G}(X_{i,-}, (P_i \setminus Y)^-)+e_{G}(X_{i,+}, (P_{i}\setminus Y)^+) + O(|X_{i,-}|+|X_{i,+}|) \Big).\]
By definition,  $|P_{i}|=|P|/k=O(\frac{n}{\log n}(\log\log n)^{1/2})$
and $|X_{i,-}|=O(n)$, $|X_{i,+}|=O(n)$. Thus, we can use Proposition \ref{prop:edges} to get
\begin{align*}
|X| \cdot |V\setminus P|\cdot p & +o(n^2p)+ \\
& \sum_{i=1}^{k}\left(|X_{i,-}|
\cdot|(P_{i}\setminus Y)^-|\cdot p+|X_{i,+}|\cdot |(P_{i}\setminus Y)^+|\cdot p+o\left(\frac{n^{2}p}{\log n}\cdot(\log\log n)^{1/2}\right)\right).
\end{align*}
Since $X=\hat{X_{i}}\cup X_{i,+}\cup X_{i,-}$ and $\left| |P_i \setminus Y| - |(P_i \setminus Y)^{-}| \right| \le 1$ (also
for $(P_i \setminus Y)^+$), this equals to
\[ |X|\cdot |V\setminus P|\cdot p+\Big(\sum_{i=1}^{k}|X\setminus\hat{X_{i}}|\cdot|P_{i}\setminus Y|\cdot p\Big) +o(n^2p).\]
As observed above, $|X\setminus\hat{X_{i}}| = (1-o(1))|X|$ for all but
$o(k)$ of indices $i$, and hence this expression becomes
$$
|X|\cdot |V\setminus P|\cdot p+
|X|p \cdot \sum_{i=1}^{k} |P_{i}\setminus Y| -o(k) \cdot \frac{|P|}{k}|X|p+ o(n^2p)
= |X| \cdot |V\setminus Y| \cdot p + o(n^2p).$$
On the other hand, this is at most the number of edges incident with $X$ in the graph $H$ which we removed, so it must be less than $|X|\cdot (\frac{1}{2}-2\varepsilon)np$.
Since $|X| \geq n/4$, we must have $|V\setminus Y| \leq \big(\frac{1}{2}-2\varepsilon+o(1)\big)n$ and
therefore $|Y|\ge(\frac{1}{2}+\varepsilon)n$.

\medskip
\noindent \textbf{Step 3} : {\em Rotating the other endpoint.}
\smallskip

In Steps 1 and 2, we constructed a set $S_P$ of size $|S_P| \ge (\frac{1}{2} + \varepsilon)n$ such that for all $v \in S_P$, there exists a
path of length $\ell$ which starts at $v$ and ends at $v_\ell$. For
each of these paths, we do the same process as in
Steps 1 and 2, now keeping $v$ fixed and rotating the other endpoint $v_\ell$. In
this way we can construct the sets $T_v$ required for the property
$\mathcal{RE}(\frac{1}{2}+\varepsilon)$.
\end{proof}

\subsection{Extension}

In the previous subsection, we showed that random graphs contain
subgraphs which can be used to perform the role of rotations. In
this subsection, we show that there exist subgraphs which can
perform the role of extensions.
\begin{defn}
\label{def:property_complement}
\label{def33}
Let $\delta$ be a positive
constant and let $G_{1}$ be a graph on $n$ vertices with property $\mathcal{RE}(\delta)$.
We say that a graph $G_2$ \emph{complements} $G_{1}$,
if for every path $P$ over the same vertex set as $G_1$
($P$ is not necessarily a subgraph of $G_1$),
either there exists a path longer than $P$ in $G_1\cup P$, or
there exist vertices $v\in S_{P}$ and $w\in T_{v}$
such that $\{v,w\}$ is an edge of $G_{1}\cup G_{2}$ (the sets
$S_P$ and $T_v$ are defined as in Definition \ref{def:property_rotext}).
\end{defn}

\begin{prop}
\label{prop:rotationextension}Let $\delta$ be a fixed positive
constant. For every $G_{1}\in\mathcal{RE}(\delta)$ and $G_{2}$
complementing $G_{1}$, the union $G_{1}\cup G_{2}$ is Hamiltonian.
\end{prop}
\begin{proof}
Let $P$ be the longest path in $G_{1}\cup G_{2}$.
By the definition of $\mathcal{RE}(\delta)$, there exists a
set $S_{P}$ such that for all $v\in S_{P}$, there exists a set $T_{v}$
such that for all $w\in T_{v}$, there exists a path of the same length as $P$
which starts at $v$ and ends at $w$. By the definition of $G_{2}$,
there exists $v\in S_{P}$ and $w\in T_{v}$ such that $\{v,w\}$
is an edge of $G_1 \cup G_{2}$. Therefore we have a cycle of length $|P|$
in $G_{1}\cup G_{2}$. Either this cycle is a Hamilton cycle or it is
disconnected to the rest of the graph, as otherwise it contradicts
the assumption that $P$ is the longest path. However, the latter
cannot happen since the
graph $G_1$ is connected by the definition of $\mathcal{RE}(\delta)$.
Thus we can conclude that the cycle we found is indeed a Hamilton cycle.
\end{proof}

The next lemma is the main lemma of this subsection and says
that the random graph complements all of its subgraphs
with small number of edges.
\begin{lem}
\label{lem:complement}For every fixed positive $\varepsilon$, there
exist constants $\delta = \delta(\varepsilon)$ and
$C=C(\varepsilon)$ such that $G=G(n,p)$ a.a.s. has the following
property: for every graph $H$ of maximum degree at most
$(\frac{1}{2} - \varepsilon)np$, the graph $G' = G-H$ complements
all graphs $R \subseteq G$ which satisfy $\mathcal{RE}(\frac{1}{2} +
\varepsilon)$ and have at most  $\delta n^{2}p$ edges.
%Every subgraph $G'=G-H$, such that the maximum degree of $H$ is at most $(1/2-\varepsilon)np$, complements all subgraphs
%$R \subseteq G(n,p)$ which satisfy $\mathcal{RE}(\frac{1}{2} + \varepsilon)$ and have at most  $\delta n^{2}p$ edges.
\end{lem}
\begin{proof}
Let $\mathcal{G}$ be the family of all subgraphs of $G$ obtained
by removing at most $(\frac{1}{2}-\varepsilon)np$ edges incident to each vertex.
%Let $G'$ be some subgraph of $G$ obtained by removing
%at most $(\frac{1}{2}-\varepsilon)np$ edges incident to each vertex.
The probability that the assertion of the lemma fails is
\begin{eqnarray}
\label{eq:eq1}
\mathbb{P} &=& \mathbb{P}\Big(\bigcup_{R\in\mathcal{RE}(\frac{1}{2} + \varepsilon), |E(R)| \leq \delta n^{2}p}\Big(\{R\subseteq G\}\cap\{\textrm{some $G' \in \mathcal{G}$ does not complement
$R$}\}\Big)\Big)\\
&\leq&  \sum_{R\in\mathcal{RE}(\frac{1}{2} + \varepsilon), |E(R)| \leq \delta n^{2}p}\mathbb{P}\Big(\textrm{some $G' \in \mathcal{G}$ does not complement $R$}\,|\, R\subseteq G\Big)\cdot
\mathbb{P}(R\subseteq G) \nonumber,
\end{eqnarray}
where the union (and sum) is taken over all labeled graphs $R$
on $n$ vertices which has property $\mathcal{RE}(\frac{1}{2} + \varepsilon)$ and at most $\delta n^{2}p$ edges.

Let us first examine the term $\mathbb{P}\Big(\textrm{some $G' \in
\mathcal{G}$ does not complement $R$}\,|\, R\subseteq G\Big)$. Let
$R$ be a fixed graph with property $\mathcal{RE}(\frac{1}{2} +
\varepsilon)$, and $P$ be a fixed path on the same vertex set. The
number of such paths is at most $n \cdot n!$, since there are $n$
choices for the length of path $P$ and there are at most
$n(n-1)\dots(n-i+1)$ paths of length $i, 1 \leq i \leq n$. If in
$R\cup P$ there is a path longer than $P$, then the condition of
Definition \ref{def33} is already satisfied. Therefore we can assume
that there is no such path in $R\cup P$. Then, by the definition of
property $\mathcal{RE}(\frac{1}{2} + \varepsilon)$, we can find a
set $S_{P}$ and for every $v \in S_P$ a corresponding set $T_{v}$,
both of size $\big(\frac{1}{2}+\varepsilon\big)n$, such that for
every $w\in T_{v}$, there exists a path of the same length as $P$ in
$R\cup P$ which starts at $v$ and ends at $w$. If there exists a
vertex $v\in S_{P}$ and $w\in T_{v}$ such that $\{v,w\}$ is an edge
of $R$, then this edge is also in $R \cup G'$ (for every $G'\in
\mathcal{G}$) and again Definition \ref{def33} is satisfied. If
there are no such edges of $R$, then since $R$ is a labeled graph,
conditioned on $R\subseteq G$, each such pair of vertices is an edge
in $G$ independently with probability $p$. Let $S'_P$ be an
arbitrary subset of $S_P$ of size $\frac{\varepsilon}{2}n$, and for
each $v \in S'_P$, define $T'_v$ to be the set $T_v \setminus S'_P$.
Since $|T'_{v}|\ge(\frac{1}{2}+\frac{\varepsilon}{2})n$, by
Chernoff's inequality, for a fixed vertex $v\in S'_{P}$, the
probability that in $G(n,p)$ this vertex has less than
$\frac{1}{2}np$ neighbors in $T'_{v}$ is at most
$e^{-\Omega_{\varepsilon}(np)}$. Since $S'_P$ is disjoint from all
the sets $T'_v$, these events are independent for different
vertices. Thus, using that $|S'_{P}| = \frac{\varepsilon}{2}n$, we
can see that the probability that all vertices $v\in S'_{P}$ have
less than $\frac{1}{2}np$ neighbors in $T'_{v}$ is at most
$e^{-\Omega_{\varepsilon}(n^{2}p)}$.

Note that if some vertex $v \in S'_P$ has at least $\frac{1}{2}np$
neighbors in $T'_{v}$, then since every $G' \in \mathcal{G}$ is
obtained from $G$ by removing at most $(\frac{1}{2}-\varepsilon)np$
edges from each vertex, there must exist a vertex $w\in T'_{v}$ such
that $\{v,w\}$ is an edge in $G'$. Therefore if some $G' \in
\mathcal{G}$ does not complement the graph $R$, then a.a.s.~there
exists some path $P$ such that all  vertices $v\in S'_{P}$ have less
than $\frac{1}{2}np$ neighbors in $T'_{v}$. Taking the union bound
over all choices of path $P$, we see that for large enough
$C=C(\varepsilon)$ and $p\ge\frac{C\log n}{n}$
\[
\mathbb{P}\Big(\textrm{some $G'$ does not complement}\: R\,|\, R\subseteq G\Big)\le n\cdot n!\cdot e^{-\Omega_{\varepsilon}(n^{2}p)}=e^{-\Omega_{\varepsilon}(n^{2}p)}.
\]
Therefore in \eqref{eq:eq1}, the right hand side can be bounded by
\[ \mathbb{P} \leq e^{-\Omega_{\varepsilon}(n^{2}p)}\cdot\sum_{R\in\mathcal{RE}(\frac{1}{2} + \varepsilon), |E(R)| \leq \delta n^{2}p}\mathbb{P}(R\subseteq G).\]
Also note that for a fixed labeled graph $R$ with $k$ edges $\mathbb{P}( R \subseteq G(n,p))=p^k$. Therefore,
by taking the sum over all possible graphs $R$ with at most $\delta n^{2}p$
edges, we can bound the probability that the assertion of the lemma fails by \[ \mathbb{P} \leq
e^{-\Omega_{\varepsilon}(n^{2}p)} \sum_{k=1}^{\delta n^{2}p}{{n \choose 2} \choose k}p^{k}
\le e^{-\Omega_{\varepsilon}(n^{2}p)} \sum_{k=1}^{\delta n^{2}p}\Big(\frac{en^{2}p}{k}\Big)^{k}.\]
For $\delta \le 1$, the summand is monotone increasing in the range $1 \le k \le \delta n^2 p$,
and thus we can take the case $k=\delta n^{2}p$ for an upper bound on every term.
This gives
\[
\mathbb{P} \le e^{-\Omega_{\varepsilon}(n^{2}p)}\cdot(\delta n^{2}p)\cdot\Big(e\delta^{-1}\Big)^{\delta n^{2}p}=e^{-\Omega_{\varepsilon}(n^{2}p)}e^{O(\delta\log(1/\delta)n^{2}p)},
\]
which is $o(1)$ for sufficiently small $\delta$ depending on $\varepsilon$. This completes the
proof.
\end{proof}

\section{Proof of the main theorem}

In this section we prove the main theorem. In view of
Lemmas \ref{lem:rotation} and \ref{lem:complement}, we can find
both the graphs we need to perform rotations and extensions.
However, we cannot immediately apply the two lemmas together,
since in order to have valid `extensions' in Lemma \ref{lem:complement},
we need the `rotation graph' to have at most $\delta n^2p$ edges.
Thus to complete the proof, we find a `rotation graph' which
has at most $\delta n^2 p$ edges.

\begin{lem}
\label{lem:splitgraph}For every positive $\varepsilon$ and $\delta<1$, there exists a constant $C=C(\varepsilon,\delta)$
such that for $p\ge \frac{C \log n}{n}$, the random graph $G(n,p)$ a.a.s. has
the following property.
For every graph $H$ of maximum degree at most $(\frac{1}{2} - 3\varepsilon)np$,
the graph $G' = G(n,p) -H$ contains a subgraph with at most $\delta n^2p$ edges satisfying $\mathcal{RE}(\frac{1}{2} + \varepsilon)$.
%Every $G'=G(n,p)-H$, such that the maximum degree of $H$ is at most $(1/2-3\varepsilon)np$,
%contains a subgraph with at most $\delta n^{2}p$ edges satisfying $\mathcal{RE}(\frac{1}{2} + \varepsilon)$.
\end{lem}
\begin{proof}
Let $C'$ be a sufficiently large constant such that for $p\ge
\frac{C' \log n}{n}$, the assertions of Proposition
\ref{prop:degreeandedges} and Lemma \ref{lem:rotation} a.a.s.~hold,
and let $C=C'/\delta$. Let $p'=\delta p$ and let $\hat{G}$ be the
graph obtained from $G(n,p)$ by taking every edge of $G$
independently with probability $\delta$. We want to analyze two
properties of $\hat{G}$ which together will imply our claim.

Call $\hat{G}$ \emph{good} if it has at most $n^2p' = \delta n^{2}p$
edges, and all of its subgraphs obtained by removing at most
$(\frac{1}{2}-2\varepsilon)np'$ edges incident to each vertex
satisfy $\mathcal{RE}(\frac{1}{2} + \varepsilon)$. Otherwise call it
{\em bad}. Note that by definition, the edge distribution of
$\hat{G}$ is identical to that of $G(n,p')$, and therefore by
Proposition \ref{prop:degreeandedges} and Lemma \ref{lem:rotation},
the probability that $\hat{G}$ is good is $1-o(1)$.  Let
$\mathcal{P}$ be the collection of graphs $G$ for which
$\mathbb{P}(\hat{G}\:\textrm{is good}\:|\: G(n,p)=G)\ge\frac{3}{4}$.
Since\begin{eqnarray*} o(1)=\mathbb{P}(\hat{G} \:\textrm{is bad}) &
\ge &
\mathbb{P}(G(n,p)\notin\mathcal{P})\cdot\mathbb{P}(\hat{G}\,\textrm{is
bad}\:|\:
G(n,p)\notin\mathcal{P})\ge\frac{1}{4}\mathbb{P}(G(n,p)\notin\mathcal{P}),\end{eqnarray*}
we know that $\mathbb{P}(G(n,p)\notin\mathcal{P})=o(1)$, or in other
words, $\mathbb{P}(G(n,p)\in\mathcal{P})=1-o(1)$. Thus from now on,
we condition on the event that $G(n,p)\in\mathcal{P}$.

Let $H$ be a graph over the same vertex set as $G(n,p)$ which has maximum
degree at most $(\frac{1}{2}-3\varepsilon)np$. Using the concentration
of hypergeometric distribution and taking union bound over all vertices of $H$, we have that with probability $1-o(1)$ the graph
$\hat{G}\cap H$ has maximum degree at most $(\frac{1}{2}-2\varepsilon)np'$.

For an arbitrary choice of $H$, since
$\hat{G}$ is good with probability at least $\frac{3}{4}$, and
$\hat{G}\cap H$ has maximum degree at most $(\frac{1}{2}-2\varepsilon)np'$
with probability $1-o(1)$,
there exists a choice of $\hat{G}$ which
satisfies these two properties. For such $\hat{G}$, by the definition of good, the graph
$\hat{G}-H$ satisfies $\mathcal{RE}(\frac{1}{2} + \varepsilon)$.
Moreover, $\hat{G}$ has at most $\delta n^{2}p$ edges and hence so does $\hat{G}-H$.
Since $\hat{G}-H \subseteq G(n,p)-H$, this proves the claim.
\end{proof}

The main result of the paper easily follows from the facts we have established so far.

\begin{proof}[Proof of Theorem \ref{thm:mainthm}]
Let $\delta$ be sufficiently small and $C$ be sufficiently large
constants such that the random graph $G=G(n,p)$ with $p\ge
\frac{C\log n}{n}$ a.a.s. satisfies Proposition
\ref{prop:degreeandedges} with $\varepsilon/2$ instead of
$\varepsilon$, and the assertions of Lemmas \ref{lem:complement} and
\ref{lem:splitgraph} with $\varepsilon/6$ instead of $\varepsilon$.
Condition on these events.

By Proposition \ref{prop:degreeandedges}, $G(n,p)$ has maximum
degree at most $(1+\frac{\varepsilon}{2})np$, and thus every subgraph of
$G(n,p)$ of minimum degree at least $(\frac{1}{2} + \varepsilon)np$ can be
obtained by removing a graph $H$ of maximum degree at most $(\frac{1}{2} - \frac{\varepsilon}{2})np$.
Thus it suffices to show that for every graph
$H$ on $n$ vertices with maximum degree at most
$(\frac{1}{2}-\frac{\varepsilon}{2})np$, the graph $G(n,p)-H$ is Hamiltonian.

Let $H$ be a graph as above.
By Lemma \ref{lem:splitgraph},
there exists a subgraph of $G(n,p)-H$ which has at most $\delta n^{2}p$
edges and has property $\mathcal{RE}(\frac{1}{2} + \frac{\varepsilon}{6})$. By
Lemma \ref{lem:complement}, $G(n,p)-H$ complements this subgraph. Therefore,
by Proposition \ref{prop:rotationextension}, $G(n,p)-H$ is Hamiltonian.
\end{proof}

\section{Concluding remarks}

In this paper, we proved that when $p \gg \log n /n$, every subgraph of the random graph $G(n,p)$ with minimum degree at least
$(1/2 +o(1))np$ is Hamiltonian. This shows that $G(n,p)$ has local resilience $(1/2 +o(1))np$ with respect to Hamiltonicity
and positively answers the question of Sudakov and Vu. It would be very interesting to better understand  the resilience of
random graphs  for values of edge probability more close to $\log n /n$, which is a threshold for Hamiltonicity.
To formalize this question we need some definitions from \cite{BeKrSu2}.

Let ${\bf a}=(a_1, \ldots, a_n)$ and ${\bf b}=(b_1, \ldots, b_n)$ be two sequences of $n$ numbers. We write $\bf{a} \leq \bf{b}$ if
$a_i \leq b_i$ for every $1 \leq i \leq n$. Given a labeled graph $G$ on $n$ vertices we denote its degree sequence by ${\bf d}_G=(d_1, \ldots, d_n)$.

\begin{defn} Let $G=([n],E)$ be a graph. Given a sequence ${\bf k}=(k_1,\ldots,k_n)$ and a monotone increasing graph property $\cal{P}$, we say that $G$ is $\bf{k}$-resilient with respect to
the property $\cal P$ if for every subgraph $H\subseteq G$ such that ${\bf d}_H\leq{\bf k}$, we have $G-H\in {\cal P}$.
\end{defn}

\noindent
It is an intriguing open problem to get a good characterization of
sequences $\bf k$ such the random graph  $G(n,p)$ with $p$ close to $\log n /n$ is
$\bf k$-resilient with respect to Hamiltonicity. Some results in this direction were obtained in \cite{BeKrSu2}.

\vspace{0.2cm} \noindent {\bf Acknowledgment.} We would like to
thank Michael Krivelevich for inspiring discussions and
conversations. We would also like to thank Alan Frieze and the two
anonymous referees for their valuable remarks.

\bibliographystyle{plain}

\end{document}